\documentclass{article}
\pdfoutput=1
\usepackage{theorem,graphicx,amssymb,amsmath}
\usepackage[retainorgcmds]{IEEEtrantools}
\usepackage{lineno}
\usepackage[hidelinks]{hyperref}
\usepackage{breakurl}
\usepackage{color}

\theorembodyfont{\slshape}

\newtheorem{theorem}{Theorem}

\newtheorem{cor}[theorem]{Corollary}

\newtheorem{conj}{Conjecture}

\def\QED{\ensuremath{{\square}}}

\def\markatright#1{\leavevmode\unskip\nobreak\quad\hspace*{\fill}{#1}}

\newcommand\blfootnote[1]{%
  \begingroup
  \renewcommand\thefootnote{}\footnote{#1}%
  \addtocounter{footnote}{-1}%
  \endgroup
}

\newcommand{\crn}{\operatorname{cr}}
\newcommand{\crs}{\overline{\operatorname{cr}}}
\newcommand{\crp}{\widetilde{\operatorname{cr}}}

\usepackage{wasysym}
\title{An Ongoing Project to Improve the Rectilinear and the Pseudolinear Crossing Constants}

\author{Oswin Aichholzer\thanks{Institute for Software Technology, Graz University of Technology, Graz, Austria. \tt{oaich@ist.tugraz.at}}
  \and Frank Duque\thanks{Instituto de Matem\'aticas, Universidad de Antioquia, Colombia.\tt{rodrigo.duque@udea.edu.co}} 
		  \thanks{Partially supported by Conacyt of Mexico grant 253261}
\and Ruy Fabila-Monroy\thanks{Departamento de Matem\'aticas, CINVESTAV, Mexico.} \thanks{ \tt{ruyfabila@math.cinvestav.edu.mx}} \footnotemark[3] 
  \and Oscar E. Garc\'ia-Quintero \footnotemark[2] \thanks{\tt{opraupe0017@gmail.com}} \footnotemark[3] 
\and Carlos Hidalgo-Toscano\thanks{INFOTEC Centro de Investigaci\'on e Innovaci\'on en Tecnolog\'ias de la Informaci\'on y Comunicaci\'on, Mexico. \tt{carlos.hidalgo@infotec.mx }} \ \footnotemark[3]  }

\begin{document}
\date{\today}
\maketitle

\begin{abstract}
A drawing of a graph in the plane is {\it pseudolinear} if the edges of the drawing can be extended to doubly-infinite curves that form an arrangement of pseudolines, that is, any pair these curves crosses precisely once.
A special case is {\it rectilinear} drawings where the edges of the graph are drawn as straight line segments. The rectilinear (pseudolinear) crossing number
of a graph is the minimum number of pairs of edges of the graph that
cross in any of its rectilinear (pseudolinear) drawings. 
In this paper we describe
an ongoing project to continuously obtain better asymptotic upper
bounds on the rectilinear and pseudolinear crossing number of the complete graph~$K_n$. 
\blfootnote{\begin{minipage}[l]{0.18\textwidth} \vspace{-6pt}\includegraphics[trim=10cm 6cm 10cm 5cm,clip,scale=0.11]{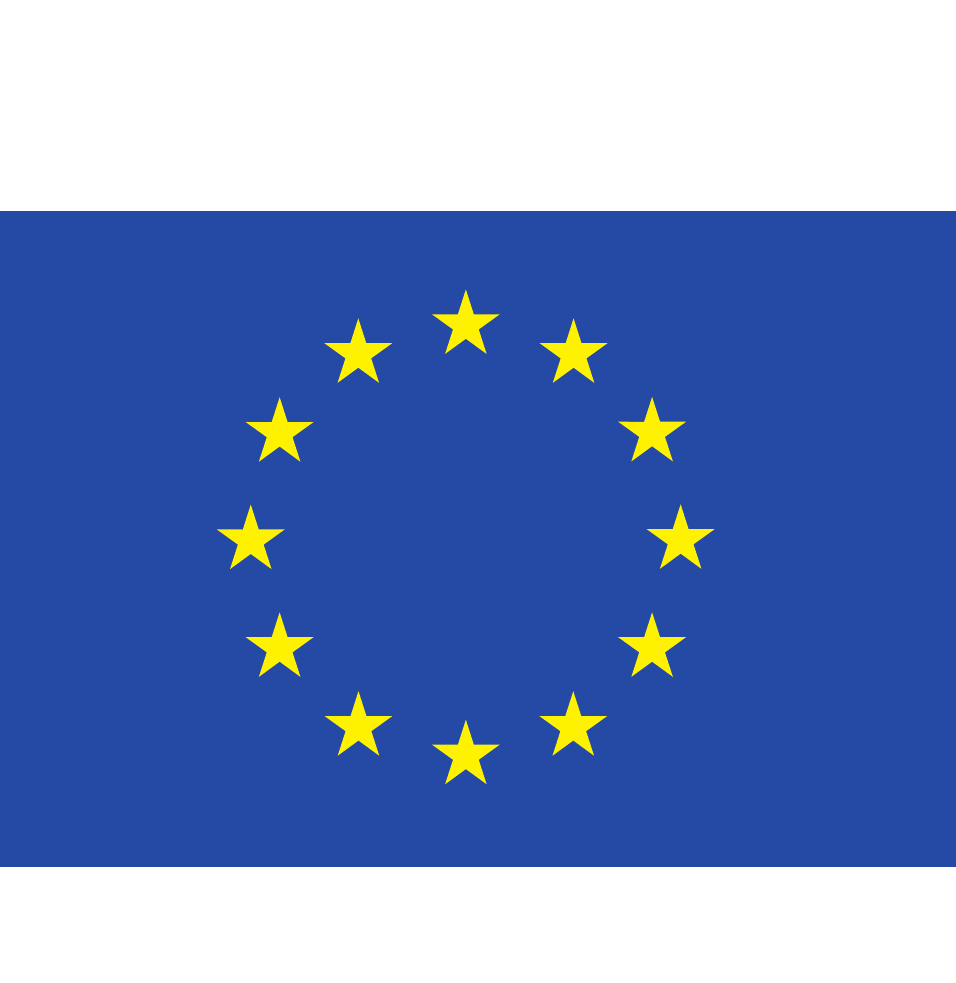} \end{minipage}  \hspace{-1.2cm} \begin{minipage}[l][1.2cm]{0.84\textwidth}\vspace{5pt}
      This project has received funding from the European Union's Horizon 2020 research and innovation programme under the Marie Sk\l{}odowska-Curie grant agreement No 734922.
     \end{minipage}}
 \end{abstract}

\section{Introduction}
Let $G$ be a graph. 
In a drawing $D$ of $G$ in the plane the vertices of $G$ are distinct points in the plane and every edge of $G$ is represented by a 
simple closed polygonal chain connecting the two points representing its vertices.
We require that in $D$ the edges do not pass through a point representing a vertex (other than their endpoints). 
A crossing in $D$ is a proper intersection of the interior of two line segments in a pair of polygonal chains
representing two edges.
Obtaining a drawing of a given graph $G$ that minimizes the number of 
crossings among all possible drawings of $G$ is a central topic in
many areas like graph drawing and discrete and computational geometry. Depending on how 
drawings are defined and on the precise way crossings are counted, there exist a huge number of
variants of crossing number problems; see~\cite{survey_marcus} for a comprehensive survey. In this work we concentrate on rectilinear and pseudolinear drawings of the complete graph~$K_n$.

The \emph{crossing number} of $G$ is the minimum number, $\crn(G)$, of pairs of edges of $G$ that cross in any drawing of $G$.
A set of points in the plane is in \emph{general position} if no three of its points are collinear.
A \emph{rectilinear drawing} of a graph $G$ is a drawing of $G$ in the plane where its vertices (points) are 
in general position and its edges are drawn as straight line segments. The \emph{rectilinear crossing number} of $G$ is the minimum 
number, $\crs(G)$,  of pairs of edges of $G$ that cross in any rectilinear drawing of $G$. 

A generalization of rectilinear drawings are pseudolinear drawings. A drawing of a graph in the plane is {\it pseudolinear} if the edges of the drawing can be extended to doubly-infinite curves that form an arrangement of pseudolines, 
that is, any pair of these curves crosses precisely once.
Similar as before, the \emph{pseudolinear crossing number}
of a graph is the minimum number, $\crp(G)$,  of pairs of edges of $G$ that
cross in any pseudolinear drawing of $G$. Clearly, $\crn(G) \le \crp(G) \le \crs(G)$. 

Finding the crossing number, rectilinear crossing number  or pseudolinear crossing number of the complete graph $K_n$ are 
important open problems in discrete geometry. 
For the crossing number of $K_n$  Harary and Hill conjectured the following. 
\begin{conj}[Harary-Hill~\cite{harary_hill}]
\[\crn(K_n)=\frac{1}{4}\left \lfloor \frac{n}{2} \right \rfloor \left \lfloor \frac{n-1}{2} \right \rfloor \left \lfloor \frac{n-2}{2} \right \rfloor \left \lfloor \frac{n-3}{2} \right \rfloor .\]
\end{conj}
According to~\cite{harary_hill} drawings achieving this bound had been found independently and Guy published them for the first time in~\cite{guy}.
The lower bounds on $\crs(K_n)$ and $\crp(K_n)$, and the upper bound on $\crn(K_n)$ together imply that for sufficiently large $n$ we have that
$\crs(K_n) \ge \crp(K_n) > \crn(K_n)$. However, it is unknown whether  $\crs(K_n) > \crp(K_n) $ (for some sufficiently large $n$).
See \cite{survey} for a nice survey.
When bounding $\crs(K_n)$ and $\crp(K_n)$ these numbers are often considered together, see e.g.~\cite{lower}; since $\crp(K_n) \leq \crs(K_n)$,
 an upper bound on $\crs(K_n)$ is also an upper bound on $\crp(K_n)$, and similar for lower bounds.

It is known that \[\lim_{n \rightarrow \infty} \frac{\crs(K_n)}{\binom{n}{4}}=q^*,\]
where  $q^*$ is a positive constant known as the \emph{rectilinear crossing constant.}
Similarly, we have \[\lim_{n \rightarrow \infty} \frac{\crp(K_n)}{\binom{n}{4}}=\widetilde{q}^*,\]
where  $\widetilde{q}^*$ is a positive constant, called the \emph{pseudolinear crossing constant.}
As it is unknown whether $\crs(K_n) > \crp(K_n)$ for sufficiently large $n$, it 
is a challenging open problem whether $q^* > \widetilde{q}^*$.

The currently best lower bounds \cite{lower} are 
\[ 0.379972 < \widetilde{q}^* \leq q^* \] 

and the previously best upper bounds \cite{pseudolinear,crossing} are
\[ \widetilde{q}^* < 0.380448 \hspace{1cm} q^*  < 0.380473.\]

In this note we describe an ongoing project to improve these upper bounds on $q^*$ and $\widetilde{q}^*$. The best upper bounds we have obtained so far\footnote{As of 2020-07-15.} are
\[q^* \le \frac{43317371729896}{113858494707069} < 0.3804491869\]
and
\[\widetilde{q}^* \le \frac{5995534434121}{15759524733750} < 0.3804387846.\]
The arXiv version of this paper [arXiv:1907.07796] will be regularly updated with new upper bounds.

\subsection{Constructing Good Drawings}

To derive upper bounds on $q^*$ one way is to produce rectilinear drawings 
of $K_n$ with few crossings for arbitrarily large values of $n$. 
The first general construction for such sets was given by Jensen~\cite{jensen} who gave explicit coordinates for the points of $K_n$. 
Around the same time Singer~\cite{singer} proposed another approach. 
His construction takes a drawing of $K_{n}$ and produces a drawing of $K_{3n}$.
If the former drawing has few crossings, then so does the latter.
Using this drawing of $K_{3n}$ and repeating the process gives a good drawing of $K_{9n}$, and so on.
This approach of iteratively generating larger sets has been successful in improving the upper bound on $q^*$ several times, see ~\cite{towards,oswinupper,bersil,upper}.
The current best iterative construction is that of \'Abrego, Cetina Fern\'andez-Merchant, Lea\~nos and Salazar(\cite{bersil} and \cite{upper}). 

A feature of these iterative constructions is that to improve the upper bound
on $q^*$ it is sufficient to find for a specific, constant value of $n$, a sufficiently
good rectilinear drawing of $K_n$. Fabila and L\'opez~\cite{crossing} used an heuristic to improve the best known rectilinear
drawing of $K_{75}$, obtaining a rectilinear drawing of $K_{75}$ with $450492$
crossings. This rectilinear drawing together with the construction of~\cite{bersil, upper} provides the until now best upper bound
of $q^* < 0.380473$.

Up to and including the work of Fabila and L\'opez~\cite{crossing} upper bounds for $\widetilde{q}^*$ were simply derived from the upper bounds of $q^*$, that is, $\widetilde{q}^* \leq q^* < 0.380473$.
Only recently Balko and  Kyn{\v{c}}l~\cite{pseudolinear} obtained the first upper bound on $\widetilde{q}^*$ which is below the upper bound of $q^*$ by showing $\widetilde{q}^* < 0.380448$.
The method they used is a generalization of~\cite{crossing} to pseudolinear drawings. They give a nice presentation of pseudolinear drawings 
with $n$-signatures\footnote{An $n$-signature is a function which assigns to every triple of vertices an orientation $\{+,-\}$.}
and show that the construction of~\cite{bersil, upper} can be adopted to the pseudolinear setting; see below for details.

Our goal is to be able to improve the upper bounds on the crossing constants in a semi-automatic way. We thus implemented the construction 
of~\cite{bersil,upper} and also its extension to pseudolinear drawings~\cite{pseudolinear}, as well as various heuristics to improve a given rectilinear or pseudolinear
drawing of $K_n$. In this paper we describe our approach in detail. In Section~\ref{sec:cons}
we briefly recall the constructions of~\cite{bersil} and ~\cite{upper} and the modification by~\cite{pseudolinear}. In Section~\ref{sec:heuristics}
we describe several heuristics we used to improve known rectilinear and pseudolinear drawings
of $K_n$. In~Section~\ref{sec:join} we describe how these tools
play together in an iterative way to regularly obtain better upper bounds on~$q^*$ and~$\widetilde{q}^*$. 

\subsection{The Construction of~\cite{bersil},~\cite{upper}, and~\cite{pseudolinear}}\label{sec:cons}

Note that a rectilinear drawing of $K_n$ is determined by the 
position of its vertices. Let $S:=\{p_1,\dots,p_n\}$ be a set 
of $n$ points in general position in the plane. Let $\crs(S)$ be the number of crossings
in the rectilinear drawing of $K_n$ where the vertices are placed at $S$. 

A \emph{halving line} of $S$ is a straight line $\ell$ passing through at 
least one point of $S$ such that
in the two open half-planes defined by $\ell$ there are the same
number of points of $S$. Note that if $n$ is odd then $\ell$ passes
only through a single point of $S$, and if $n$ is even then $\ell$ passes
through two points of $S$. 

Let $G$ be the bipartite graph with vertex partition $(A,B)$, where $A:=S$
and $B$ is the set of the halving lines of $S$. A pair
$(p,\ell)$ in $(A,B)$ is adjacent in $G$ if and only if $\ell$ passes through 
$p$. A \emph{halving matching} of $S$ is a matching of $G$, in which
every point in $A$ is matched to a halving line in $B$.
If $n$ is odd then a halving matching of $S$ always exists. In this case
every halving line contains exactly one point of $S$ and every point of $S$
is contained in at least one halving line (actually an infinite number of them). 
So any choice of a halving line for every point of $A$ is a halving matching.
A halving matching may
not exist if $n$ is even. For example, a set of $4$ points in which exactly three
of them are in the convex hull does not contain a halving matching.

A rough description of the construction of~\cite{bersil} and \cite{upper} is as follows.
Let $M=\{(p_1,\ell_1),\dots, (p,\ell_n)\}$ be a halving matching for $S$.
For every $\ell_i$ assume that $\ell_i$ is directed, and let $\vec{v}_i$ be the 
direction vector of $\ell_i$. Let $S'$ be the point set that results by replacing
each $p_i$ in $S$ with the pair of points $p_i+\varepsilon\vec{v}_i$ and
$p_i-\varepsilon\vec{v}_i$ for an arbitrarily small (but positive) value of $\varepsilon$. 

This construction was first described in~\cite{bersil}; furthermore, they showed
that if $S$ has an even number of points then  $S'$ also has a halving matching.
This allows for the iterative construction mentioned above. If $S$ has an
odd number of points then it always has a halving matching; however, it is not obvious
that $S'$ should also have a halving matching.
In~\cite{upper} they showed that in this case $S'$ also has a halving matching.
The following theorem gives the upper bound on $q^*$  that is derived from this construction. 

\begin{theorem}[\cite{upper,bersil}]\label{thm:rectilinear_constant}
 Let $S$ be a set of $n$ points in general position in the plane, such
 that $S$ has a halving matching. Then
 \[q^* \le \frac{24 \crs(S)+3n^3-7n^2+(30/7)n}{n^4}.\]
\end{theorem}

The above theorem is derived by showing that in every doubling step from $S$ to $S'$ it holds that
$\crs(S')=16 \crs(S) + (n/2)(2n^2-7n+5)$. It is not hard to see that the concept of halving matching also translates to pseudolinear drawings.
Consequently it is shown in~\cite{pseudolinear} that the construction of doubling can be applied to pseudolinear drawings and that the halving matching property is also preserved.
Thus, a similar relation holds when doubling a pseudolinear drawing $D$ to $D'$, namely $\crp(D')=16 \crp(D) + 2n \left(\lceil n/2 \rceil^2+\lfloor n/2 \rfloor^2\right)-7n^2/2+5n/2$~\cite{pseudolinear}.
For even $n$ this is the same bound as for the rectilinear case. Although not explicitly stated in~\cite{pseudolinear} this therefore leads to the following bound for~$\widetilde{q}^*$.

\begin{cor}\label{cor:pseudolinear_constant}
 Let $D$ be a pseudolinear drawing of $K_n$, such that $D$ has a halving matching. Then
  \[\widetilde{q}^* \le \frac{24 \crp(D)+3n^3-7n^2+(30/7)n}{n^4} \textrm{ for } n \textrm{ even}; and\]
  \[\widetilde{q}^* \le \frac{24 \crp(D)+3n^3-7n^2+(81/14)n}{n^4}  \textrm{ for } n \textrm{ odd} .\]
\end{cor}

While the above two results will be directly used to derive our improved upper bounds, we are also interested in actually obtaining the good sets of doubled cardinality; see also Section~\ref{sec:join}.

For rectilinear drawings the above construction has been implemented as a Python program in~\cite{oscar} and produces point sets with integer coordinates.
The choice of using integer coordinates helps to avoid numerical issues, since as long as we use a library (or programming language) that handles arbitrarily large
integers, arithmetic precision is not an issue. We therefore used the implementation of~\cite{oscar} also for our computations.

For pseudolinear drawings no numerical issues arise, as no realization with points is needed and thus no coordinates are computed.
We follow the lines of~\cite{pseudolinear} and use a representation of a pseudolinear drawing by $n$-signatures.
This is a function which assigns to every triple of vertices an orientation $\{+,-\}$.
The advantage of this representation is that it is easy to check if a given $n$-signature can be realized as a pseudolinear drawing and that the number of crossings can be computed directly from the representation.
Moreover, the above mentioned results for the doubling construction have been derived; see again~\cite{pseudolinear} for more details.
We implemented this approach also as a Python program to generate good pseudolinear drawings of twice the size.

\section{Heuristics}\label{sec:heuristics}

We now describe various heuristics that attempt to improve a given drawing of~$K_n$.
For a rectilinear drawing the vertices are the points of the point set $S$ and $p$ is a point of $S$.
For a pseudolinear drawing $D$ all triple orientations $v_iv_jv_k$, $1 \leq i < j < k \leq n$ are given.

\subsection{Moving a Point to a New Random Location}\label{sec:moving}

In~\cite{crossing} the following simple heuristic is used for rectilinear drawings.
Choose a random point $q \notin S$ in a certain neighborhood of $p$. If 
$\crs(S\setminus\{p\} \cup \{q\}) \le \crs (S)$ then
replace $p$ with $q$ in $S$. In~\cite{crossing}, $\crs(S\setminus\{p\} \cup \{q\}) $ is computed
in $O(n^2)$ time. If after some time no improvement is found then the size of the neighborhood of~$p$
from which~$q$ is chosen, is made smaller. 
We use the amortized faster algorithm of~\cite{crossing_am}.
In that algorithm a set $Q$ of $\Theta(n)$ candidate points for the new position of $p$
is chosen. The set of values \[\{\crs(S \setminus \{p\} \cup \{q\}) :q \in Q\}\]
is computed in $O(n^2)$ time. Note that this is linear per point
in $Q$. We then consider the best point in $Q$ as a possible replacement for $p$.

For pseudolinear drawings a very similar idea is used. Obviously there are no points to be moved. But observe that moving a point continuously in the rectilinear case corresponds to changing the orientation of point triples one after another (since there are no collinear points), that is, changing the order type of the set step by step. The difficulty in the rectilinear setting is that this is a geometric process, and even if all points have integer coordinates, moving the point on an integer grid might not be sufficient to get all triple changes separately. 

To the contrary, changing triple orientations in a pseudolinear drawing given as $n$-signature is trivial.
We choose a random vertex triple $v_iv_jv_k$, $1 \leq i < j < k \leq n$, of the drawing $D$ and invert its orientation.
Then we check if this new $n$-signature is still realizable as a pseudolinear drawing $D'$.
In case $D'$ exists, we keep the new orientation if $\crp(D') \leq \crp(D)$ and set $D=D'$.
Otherwise we revert the change. Iterating this process eventually gives pseudolinear drawings with less crossings.

It is interesting to observe that our experiments show that this local optimization heuristic works in average significantly better for pseudolinear drawings than for rectilinear drawings.
The main reason might be that for pseudolinear drawings the algorithm has a combinatorial flavor, while for rectilinear drawings the precise geometry of the arrangement of the point set (e.g.\ the size of cells) plays a role.
We therefore present in the next section an approach for rectilinear drawings which avoids this bottleneck. 

\subsection{The Point Set Explorer}\label{sec:point_explorer}

Consider the line arrangement, $\mathcal{A},$ spanned by the set of lines
passing through every pair of points of $S \setminus \{ p\}$. Let $C$ be a cell
of $\mathcal{A}$ and let $q$ be a point in $C$. Note that $\crs(S\setminus \{p\}\cup \{q\})$ has 
the same value regardless of the choice of the precise location of $q$ within $C$.

\begin{figure}
	\begin{center}
		\includegraphics[width=0.7\textwidth]{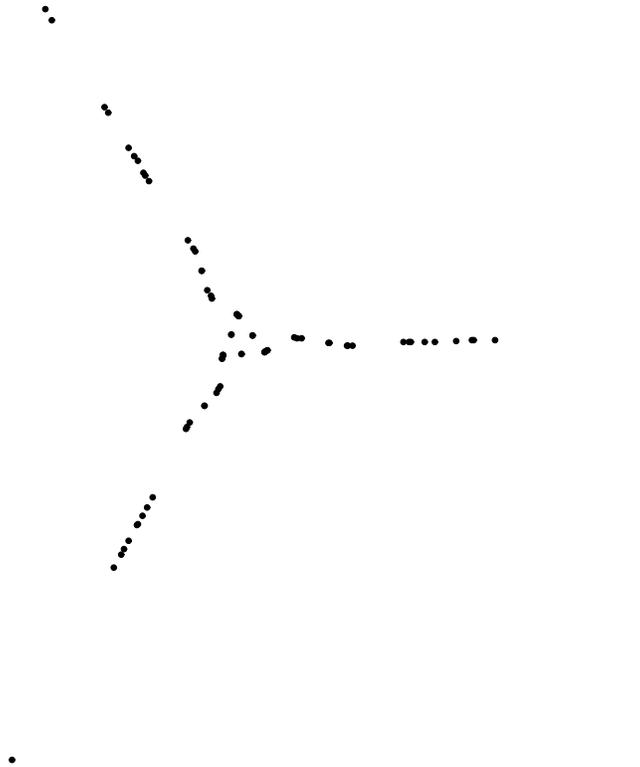}
	\end{center}
	\caption{A drawing of $K_{75}$ with $450492$ crossings}
	\label{fig:75}
\end{figure}

All the best known rectilinear drawings of $K_n$ consist of three ``arms'' with close to
$n/3$ points each; the points in each arm are close to being collinear. See Figure~\ref{fig:75} for an example.
This implies that many of the cells in the line arrangements of such sets have area close to zero. Therefore, 
in the heuristic described in~\ref{sec:moving} many of the candidate points
fall in the same cells and many cells are never visited.

In~\cite{carlos} the following algorithm is implemented. The algorithm computes the cell of $\mathcal{A}$
that contains $p$ in $O(n^2 \log n)$ time.
Afterwards, moving between adjacent cells, and finding a point in the interior of these cells, takes $O(\log^2 n)$ time. 

Using this implementation we produce a sequence $C=C_1,\dots,C_m$ of consecutive adjacent cells. At each
$C_i$ we choose a point $q_i$ as a candidate to replace $p$ in $S$.  Since we are moving
between adjacent cells, $\crs(S \setminus \{p\} \cup \{q_i\})$  can actually be computed in $O(1)$ time.

In this way the size of the cells in the arrangement and thus the precise geometry of our point set does not play a role anymore.
We obtain an algorithm with a more combinatorial behavior than with the previous heuristic.

\subsection{Finding Good Subdrawings} \label{sec:subdrawing}

It is often the case that for $m > n$, drawings of $K_m$ with few crossings contain drawings of $K_n$ with few crossings.
We observed that when optimizing drawings of $K_m$ with the heuristics presented in the previous section, they often contain better subdrawings of $K_n$ than before the optimization.
In other words, it is possible to improve drawings of $K_n$ by optimizing drawings of $K_m$.
It might sound counter intuitive to optimize the larger set.
But observe that there are many subsets of size $n$ in a set of size $m$. 
So if just one of these subdrawings is improved, this approach is successful.

Given a drawing of $K_m$, we use the following
heuristics to find drawings of $K_n$ with few crossings.

\subsubsection*{Removing one point at a time}
For rectilinear drawings we remove one point of $S$ at a time, to obtain
smaller drawings. 
We use an implementation of an algorithm described in~\cite{crossing_am}
that does the following. In $O(n^2)$ time, it computes the set of values 
\[\{\crs(S\setminus \{p\}): p \in S\}.\]
We remove the point $p$ that minimizes $\crs(S\setminus \{p\})$, thus finding
a rectilinear drawing of $K_{m-1}$. We proceed iteratively
in this way to find rectilinear drawings of $K_n$ with $n=m-1,m-2,\dots,27$ points.

For pseudolinear drawings this is done in a similar way.

\subsubsection*{Removing two or more points at a time}
It is known that point sets minimizing the rectilinear crossing number have a triangular convex hull, and most of them have several layers consisting of three vertices each.
Moreover, during the process of generating good examples to get an improved crossing constant we could observe that the cardinality of the best sets seems to be a multiple of 3.
So actually it turned out that sometimes it is better to remove more than one point at each step.
That means we look at a tuple or triple of points which, when being removed at once, reduces the crossing number the most.
Occasionally, this provides rectilinear drawings that cannot be obtained by removing one point at a time as in the previous heuristic.
However, the amortized speed up described above is lost, and thus this method needs more computational resources. 

To our own surprise good pseudolinear drawings show precisely the same behavior.
There also the best drawings (w.r.t.\ the obtainable pseudolinear crossing constant) we found so far have a cardinality which is a multiple of 3.
We therefore use the same approach of removing a tuple or triple of vertices at the same time also for pseudolinear drawings.

\section{Joining Everything Together}\label{sec:join}

Several computers run all the heuristics described in Section~\ref{sec:heuristics} permanently in the background.
Processes are started for the best drawings of $K_n$ for some specific values of $n$ - typically the currently best three or four cardinalities for both types of drawings, rectilinear and pseudolinear. The newly-found drawings are sent daily to
a central node.  If after for some time, no new improvements are made on the upper bound of either $q^*$ or $\widetilde{q}^*$ , we take the corresponding drawing of $K_n$
 providing the best upper bound and apply the doubling construction described in Section~\ref{sec:cons} to obtain a drawing of $K_{2n}$ providing the same upper bound on $q^*$ or $\widetilde{q}^*$, respectively. 
From this set we then compute good subsets by removing points/vertices as described in Section~\ref{sec:subdrawing}. In parallel we locally optimize the new drawing $K_{2n}$ and all the obtained subdrawings by the methods described in Sections~\ref{sec:moving} and~\ref{sec:point_explorer}. The steps of local optimization and computing subdrawings are interleaved and iterated. Note, however, that the local optimization is applied only for a limited time. The best sets obtained in these iterations serve as new starting sets of our optimization, and the whole process is restarted.
 
\subsection{Results for Rectilinear Drawings}

With the just described approach we were able to obtain good rectilinear drawings of $K_n$ for up to $n = 3240$ points.
The best crossing constant is obtained by using a rectilinear drawing of $K_{2643}$ which has $771218714414$ crossings.
The interested reader can download the file with the coordinates of the points 
from \href{http://www.crossingnumbers.org/projects/crossingconstants/rectilinear.php}{\small \tt http://www.crossingnumbers.org/projects/crossingconstants/rectilinear.php}.
The reason why larger sets do not necessarily provide a better constant lies in the before mentioned search for good subsets, which is an essential ingredient of our mixed heuristic.
Plugging the values of this drawing into Theorem~\ref{thm:rectilinear_constant} we obtain the following result.

\begin{theorem}
\[q^* \le \frac{43317371729896}{113858494707069} < 0.3804491869\]
\end{theorem}

\subsection{Results for Pseudolinear Drawings}

For pseudolinear drawings we have been able to get good results for $K_n$ for up to $n = 2502$.
The reason why the maximal cardinality is smaller than the one for rectilinear drawings is that here we have to store every triple orientation explicitly, which needs several hundreds MB for each set.
The best crossing constant is obtained by using a pseudolinear drawing of $K_{2205}$ which has $373382224051$ crossings.
The interested reader can download the file with the information of the triple orientation 
from \href{http://www.crossingnumbers.org/projects/crossingconstants/pseudolinear.php}{\small \tt http://www.crossingnumbers.org/ projects/crossingconstants/pseudolinear.php}.
Plugging the values of this drawing into Corollary~\ref{cor:pseudolinear_constant} we obtain the following result.

\begin{theorem}
\[\widetilde{q}^* \le \frac{5995534434121}{15759524733750} < 0.3804387846.\]
\end{theorem}

%
\end{document}